\documentclass[12pt]{amsart}
\usepackage{amssymb,latexsym}
\makeatletter
\let\atopwithdelims\@@atopwithdelims
\let\over\@@over
\newtheorem{theorem}{Theorem}[section]
\newtheorem{proposition}[theorem]{Proposition}
\newtheorem{corollary}[theorem]{Corollary}

\newtheorem{lemma}[theorem]{Lemma}

\newtheorem{conjecture}[theorem]{Conjecture}
\newcounter{enumit}
\newenvironment{enumit}{\begin{list}{({\it\roman{enumit}})}{
                        \usecounter{enumit}}
                      }{\end{list}}

\newcommand{\for}{\textrm{for}}

\begin{document}
\title[Extended Linial Arrangements for Root Systems]{Extended Linial 
Hyperplane Arrangements\\for Root Systems\\and a Conjecture of Postnikov 
and Stanley}
\author{Christos~A.~Athanasiadis}
\address{\hskip-\parindent Christos~A.~Athanasiadis\\
Department of Mathematics, University of 
Pennsylvania, 209 South 33rd Street, Philadelphia, PA 19104-6395}
\curraddr[until August 10, 1997]{\sc
Mathematical Sciences Research Institute\\
1000 Centennial Drive\\
Berkeley, CA 94720}

\email{athana@msri.org}
\date{\today}
\thanks{Supported by a postdoctoral fellowship from the Mathematical 
Sciences Research Institute, Berkeley, California. Research at MSRI 
is supported in part by NSF grant DMS-9022140.}
\begin{abstract}
A hyperplane arrangement is said to satisfy the ``Riemann hypothesis'' 
if all roots of its characteristic polynomial have 
the same real part. This property was conjectured by Postnikov and 
Stanley for certain families of arrangements which are defined for any 
irreducible root system and was proved for the root system $A_{n-1}$. 
The proof is based on an explicit formula \cite{Ath1, Ath2, PS} for 
the characteristic polynomial, which is of independent combinatorial 
significance. Here our previous derivation of this formula is simplified 
and extended to similar formulae for all but the exceptional root 
systems. The conjecture follows in these cases.
\end{abstract}

\maketitle

\section{Introduction}

Let ${\mathcal A}$ be a hyperplane arrangement in ${\mathbb R}^n$, i.e. 
a finite collection of affine subspaces of ${\mathbb R}^n$ of codimension 
one. The \textit{characteristic polynomial} \cite[\S 2.3]{OT} of 
${\mathcal A}$ is defined as
\[ \chi({\mathcal A}, q) = \sum_{x \in L_{\mathcal A}} \mu (\hat{0},x) \ 
q^{\dim x}, \]
where $L_{\mathcal A}$ is the poset of all affine subspaces of 
${\mathbb R}^n$ which can be 
written as intersections of some of the hyperplanes of $\mathcal A$, 
$\hat{0} = {\mathbb R}^n$ is the unique minimal element of $L_{\mathcal A}$ 
and $\mu$ stands for its M\"obius function \cite[\S 3.7]{St1}. The polynomial 
$\chi({\mathcal A}, q)$ is a fundamental combinatorial and topological 
invariant of $\mathcal A$ and plays a significant role throughout the 
theory of hyperplane arrangements \cite{OT}.

\smallskip
Very often the polynomial $\chi({\mathcal A}, q)$ factors completely over 
the nonnegative integers. This happens, for instance, when $\mathcal A$ is
a \textit{Coxeter arrangement}, i.e. the arrangement of reflecting 
hyperplanes of a finite Coxeter group \cite[p. 3]{OT}. A number of theories 
\cite{St0, Te, JP} have been developed to explain this phenomenon (see also
the survey \cite{Sa}). A different phenomenon has 
appeared in recent work of Postnikov and Stanley \cite{PS} and is referred
to as the ``Riemann hypothesis'' for $\mathcal A$. It asserts that all roots
of $\chi({\mathcal A}, q)$ have the same real part. This property was 
conjectured in \cite{PS} for certain affine \textit{deformations} of 
Coxeter arrangements and was proved for the Coxeter type $A_{n-1}$. The
proof was based on explicit formulae for the characteristic polynomials,
first obtained in \cite{Ath1} \cite[Part II]{Ath2}. In this paper we 
improve and extend our previous arguments to treat case by case all but 
the exceptional Coxeter types. There is no general theory known that 
could give a more uniform proof.
 
\smallskip
We first state precisely the Conjecture of Postnikov and Stanley and 
our main result.

\vspace{0.1 in}
\textit{The main result}. A \textit{root system} $\Phi$ will be a 
crystallographic root system \cite[\S 2.9]{Hu} which is not necessarily 
reduced, i.e. if $\alpha, \beta \in \Phi$ with $\beta = c \alpha$ then we
do not require that $c = \pm 1$. This includes the non-reduced system 
$BC_n$ which is the union of $B_n$ and $C_n$. Let $\mathcal A$ be the 
Coxeter arrangement corresponding to $\Phi$. A deformation of $\mathcal A$ 
\cite{PS, St2} is an arrangement each of whose hyperplanes is parallel to 
some hyperplane of $\mathcal A$. Fix a system of positive roots $\Phi^+$
and let $a \leq b$ be integers. We denote by $\hat{\mathcal A}^{[a, b]} \, 
(\Phi)$ the deformation of $\mathcal A$ which has hyperplanes 
\[ (\alpha, x) = k \ \ \for \ \ \alpha \in \Phi^+ \ \ \textrm{and} \ \ 
k = a, a+1,\ldots,b. \]
This reduces to $\mathcal A$ if $a=b=0$. The conjecture of Postnikov and 
Stanley from \cite[\S 9]{PS} is as follows.
\begin{conjecture}
Let $\Phi$ be an irreducible root system in ${\mathbb R}^l$ and $a, b$ be 
nonnegative integers, not both zero, satisfying $a \leq b$. If 
$h_{\hat{\mathcal A}}$ is the number of hyperplanes of $\hat{\mathcal A} 
= \hat{\mathcal A}^{[-a+1, b]} \, (\Phi)$ then all roots of 
$\chi(\hat{\mathcal A}, q)$ have real part equal to $h_{\hat{\mathcal A}}/ 
\, l$.
\end{conjecture}
\textit{Note}. For any arrangement ${\mathcal A}$ in ${\mathbb R}^l$,
the sum of the roots of $\chi({\mathcal A}, q)$ is equal to the number
$h_{{\mathcal A}}$ of hyperplanes of ${\mathcal A}$. Hence, if all roots 
of $\chi({\mathcal A}, q)$ have the same real part, this has to be
$h_{\mathcal A}/ \, l$.

\smallskip
The characteristic polynomial of $\hat{\mathcal A}^{[a, b]} \, (\Phi)$ is 
independent of the choice of positive roots $\Phi^+$, so from now and on 
we assume that this set is as in \cite[\S 2.10]{Hu}. We then abbreviate
$\hat{\mathcal A}^{[a, b]} \, (\Phi)$ as $\hat{\mathcal B}_n ^{[a, b]}$,
$\hat{\mathcal C}_n ^{[a, b]}$, $\hat{\mathcal D}_n ^{[a, b]}$ or
$\hat{\mathcal {BC}}_n ^{[a, b]}$ if $\Phi = B_n$, $C_n$, $D_n$ or $BC_n$
respectively. For $\Phi = A_{n-1}$, this is an arrangement in 
${\mathbb R}^{n-1}$. For convenience, we denote by
$\hat{\mathcal A}_n ^{[a, b]}$ the arrangement of hyperplanes in 
${\mathbb R}^n$ of the form 
\[ x_i - x_j = a, a+1,\ldots,b \ \ \for \ \ 1 \leq i < j \leq n, \]
so that $\hat{\mathcal A}_n ^{[a, b]}$ is the product 
\cite[Definition 2.13]{OT} of the empty one dimensional arrangement
and $\hat{\mathcal A}^{[a, b]} \, (\Phi)$ and hence 
$\chi(\hat{\mathcal A}_n ^{[a, b]}, q) = q \, \chi(\hat{\mathcal A}
^{[a, b]} \, (\Phi), q)$, where $\Phi = A_{n-1}$. The arrangements 
$\hat{\mathcal A}_n ^{[1,b]}$
are referred to as the \textit{extended Linial arrangements}. They were 
studied enumeratively because of a remarkable conjecture of Linial and 
Stanley, first
proved by Postnikov \cite[Thm. 8.2]{PS} (also in \cite[\S 4]{Ath1} 
\cite[\S 6.4]{Ath2}), about the number of regions of the \textit{Linial 
arrangement}, the one which corresponds to $b=1$. The polynomials 
$\chi(\hat{\mathcal A}_n ^{[1, b]}, q)$ were first computed explicitly 
in \cite[\S 4]{Ath1} \cite[\S 6.4]{Ath2} with the finite field method.
We use the same method to find similar explicit formulae in the case of 
the other classical root systems and prove the following theorem.  
\begin{theorem}
Conjecture 1.1 holds for the infinite families of root systems $A_{n-1}$, 
$B_n$, $C_n$, $D_n$ and $BC_n$, where $n \geq 2$.
\end{theorem}
As remarked earlier, the proof of Theorem 1.2 will be done case by case. 
No uniform proof is known.  

\vspace{0.1 in}
The paper is organized as follows: Section 2 contains a review and 
refinement of the finite field method of \cite{Ath1} \cite[Part II]{Ath2} 
and other useful background. In Section 3 we simplify substantially the 
derivations of the formulae for $\chi(\hat{\mathcal A}_n ^{[1, b]}, q)$ and 
$\chi(\hat{\mathcal A}_n ^{[0, b]}, q)$ given in \cite[\S 4]{Ath1} 
\cite[\S 6.4]{Ath2}. In particular, we get a simple proof of Postnikov's 
theorem for the number of regions of the Linial 
arrangement. In Section 4 we obtain similar formulae for the root systems
$B_n$, $C_n$, $D_n$ and $BC_n$. In Section 5 we use the results of Sections
3 and 4 and an elementary lemma, employed by Postnikov and Stanley, to 
complete the proof of Theorem 1.2. We conclude with some remarks in 
Section 6.

\section{Background}

We first review the finite field method of \cite{Ath1} \cite[Part II]{Ath2}. 
This method reduces the computation of the characteristic polynomial to 
a simple counting problem in a vector space over a finite field. It will 
be more convenient here to work over the abelian group ${\mathbb Z}_q$ 
of integers modulo $q$, where $q$ is not necessarily a power of a prime. 
We will naturally restrict our attention to hyperplane arrangements, as 
opposed to the more general \textit{subspace} arrangements \cite{Bj,BjE}. 

\smallskip
Let ${\mathcal A}$ be any hyperplane arrangement in ${\mathbb R}^n$ 
and $q$ be a positive integer. We call ${\mathcal A}$ a 
${\mathbb Z}$-\textit{arrangement} if its hyperplanes are given by 
equations with integer coefficients. Such equations define subsets 
of the finite set ${\mathbb Z}_q ^n$ if we reduce their coefficients 
modulo $q$. We denote by $V_{\mathcal A}$ the union of these subsets, 
supressing $q$ from the notation. The next theorem is a variation of 
\cite[Thm.\ 2.2]{Ath1} \cite[Thm.\ 5.2.1]{Ath2} (see also the original 
formulation in \cite[\S 16]{CR} as well as \cite[Thm.\ 2.69]{OT}, 
\cite[Thm.\ 2.1]{BS} and Proposition 3.2 and Lemma 5.1 in \cite{BjE}).
\begin{theorem}
Let ${\mathcal A}$ be a ${\mathbb Z}$-hyperplane arrangement in 
${\mathbb R}^n$. There exist positive integers $m, k$ which depend 
only on ${\mathcal A}$, such that for all $q$ relatively prime to 
$m$ with $q > k$,
\[ \chi({\mathcal A},q) = \# \left({\mathbb Z}_q ^n - V_{\mathcal A}
\right). \]
\end{theorem}
\begin{proof}
Let $H_1, H_2,\ldots,H_r$ be some of the hyperplanes of ${\mathcal A}$, 
$X \subseteq {\mathbb R}^n$ be their intersection and $X_q$ be the 
intersection of the corresponding subsets of ${\mathbb Z}_q ^n$. It 
suffices to guarantee that $\# X_q = q^{\dim X}$ if $X$ is nonempty 
and $X_q = \emptyset$ otherwise, for any such choice of hyperplanes. 
The result then follows by M\"obius inversion as in \cite{Ath1, Ath2, BS,
CR} or, equivalently, by the argument given in Propositions 3.1 and 3.2 
of \cite{BjE}. Let $X$ be described by the linear system
\begin{equation}
Ax = b,
\label{syst}
\end{equation}
where $A$ is an $r$ by $n$ ${\mathbb Z}$-matrix and $b$ has integer 
entries. Since there are invertible ${\mathbb Z}$-matrices $P, Q$ such 
that $P^{-1}AQ$ is diagonal, we can assume that (\ref{syst}) consists 
of the equations $d_i x_i = b_i$ for $1 \leq i \leq r$. It suffices to 
choose $m, k$ so that $d_i | m$ whenever $d_i \neq 0$ and $k > |b_i|$ 
whenever $d_i = 0$. 
\end{proof}

\noindent
\textit{Remark.} We can choose $m$ to be $1$ or $2$ if ${\mathcal A}$ 
is a ${\mathbb Z}$-deformation of ${\mathcal A}_n$ or ${\mathcal {BC}}_n$, 
respectively, i.e if ${\mathcal A}$ is contained in some 
$\hat{\mathcal A}_n ^{[a, b]}$ or $\hat{\mathcal {BC}}_n ^{[a, b]}$
for integers $a \leq b$. We will make use of this fact in the following 
sections without further comment. Also, we can choose $k = 0$ if 
${\mathcal A}$ is central.

\vspace{0.1 in}
\noindent
\textit{Notation.} We often write 
\[ \phi_a(y) := 1 + y + y^2 + \cdots + y^{a-1}. \]
This polynomial will appear repeatedly in the formulae of Sections 3 and 
4. The \textit{shift operator} $S$ acts on polynomials $f$ of one variable 
by 
\[ Sf (y) := f(y-1). \]

\vspace{0.1 in}
The following elementary lemma will be needed in the next sections.  
\begin{lemma}
For fixed positive integers $a, n$ let 
\[ \left(\phi_a(y)\right)^n := \left(1 + y + y^2 + \cdots + 
y^{a-1}\right)^n = \sum_{k=0}^{n(a-1)} c_k \, y^k. \]
If $0 \leq i \leq a-1$ and $f$ is a polynomial of degree less than $n$, 
then the sum
\[ \Sigma_i f(y) := \left(\sum_{k \equiv i \ ({\rm mod} \ a)} c_k \, 
S^k\right) f(y) = \sum_{k \equiv i \ ({\rm mod} \ a)} c_k \, f(y-k) \]
is independent of $i$ and hence
\[ \Sigma_i f(y) = \frac{1}{a} \, \left(\phi_a(S)\right)^n f. \]
\end{lemma}
\begin{proof}
By linearity, it suffices to prove the result for $f(y) = y^j$, where 
$0 \leq j \leq n-1$. We fix such a $j$ and $r$ with $0 \leq r \leq j$. 
The coefficient of $y^{j-r}$ in $\Sigma_i y^j$ is $(-1)^r {j \choose r} 
\, s_i$, where
\[ s_i = \sum_{k \equiv i \ ({\rm mod} \ a)} c_k \, k^r. \]
Therefore, it suffices to show that $s_0 = s_1 = \cdots = s_{a-1}$. Note 
that 
\[ \sum_{k=0}^{n(a-1)} c_k \, k^r \, y^k = \left(y \, \frac{d}{dy} 
\right)^r (\phi_a(y))^n \]
is divisible by $\phi_a(y)$. Thus, setting $y=\omega$, a primitive $a$th 
root of unity, we get $s_0 + s_1 \omega + s_2 \omega^2 + \cdots + s_{a-1} 
\omega^{a-1} = 0$. The same is true if $\omega$ is replaced with 
$\omega^m$ for $m=2,\ldots,a-1$. Hence the column vector 
$(s_0, s_1,\ldots,s_{a-1})^t$ is in the kernel of the $a-1$ by $a$ 
matrix $\Omega$ whose entry in position $(m, l)$ is equal to 
$\omega^{m(l-1)}$. The first $a-1$ columns of the matrix $\Omega$ are 
linearly independent, so it has rank $a-1$ and a one dimensional kernel. 
The kernel is clearly spanned by the column vector with all entries 
equal to $1$ so indeed, $s_0 = s_1 = \cdots = s_{a-1}$. 
\end{proof}

\section{The root system $A_{n-1}$}

In this section we consider the case of the root system $A_{n-1}$. We 
rederive the formulae for $\chi(\hat{\mathcal A}_n ^{[0,a]}, q)$ and 
$\chi(\hat{\mathcal A}_n ^{[1,a]}, q)$ along the lines of 
\cite[\S 4]{Ath1} \cite[\S 6.4]{Ath2} but use a simpler and more direct
combinatorial argument. This case will serve as a prototypical example 
of application of the finite field method, which we will adjust in the 
next section to the case of other root systems. 

\vspace{0.1 in}
For the following proof, we represent an $n$-tuple $x = 
(x_1, x_2,\ldots,x_n)$ of distinct elements of ${\mathbb Z}_q$ as a 
placement of the integers $1, 2,\ldots,n$ and $q-n$ indistinguishable 
balls along a line. Such a placement corresponds to the $n$-tuple $x$ 
for which $x_i + 1$ is the position that $i$ occupies, counting from 
the left. For example, for $q = 10$ and $n = 4$, the placement
\begin{equation}
4 \ \bigcirc \, \bigcirc \, \bigcirc \ 2 \ 3 \ \bigcirc \, \bigcirc \ 
1 \ \bigcirc 
\label{pl}
\end{equation}
corresponds to the $4$-tuple $(8, 4, 5, 0)$ of elements of 
${\mathbb Z}_{10}$. We denote by $[y^k] \, F(y)$ the coefficient of
$y^k$ in the formal power series $F(y)$.
\begin{proposition}{\rm (\cite[Thm.\ 4.4]{Ath1} \cite[Thm.\ 6.4.4]{Ath2})}
For all $a \geq 1$ and $q > an$ we have
\begin{equation}
\chi(\hat{\mathcal A}_n ^{[0,a]}, q) = q \, [y^{q-n}] \
 (1 + y + y^2 + \cdots + y^{a-1})^n \sum_{j=0}^{\infty} j^{n-1} y^{aj}.
\label{ch1}
\end{equation}
\end{proposition}
\begin{proof}
Theorem 2.1 implies that, for large positive integers $q$, 
$\chi(\hat{\mathcal A}_n ^{[0,a]}, q)$ counts the number of $n$-tuples 
$x = (x_1, x_2,\ldots,x_n) \in {\mathbb Z}_q ^n$ which satisfy
\[ x_i - x_j \neq 0, 1,\ldots,a \]
in ${\mathbb Z}_q$ for all $1 \leq i < j \leq n$. Since $x$ satisfies 
these conditions if and only if $x + m := (x_1 + m,\ldots,x_n + m)$ does 
so, we can assume that, say, $x_n = 0$ and disregard the factor of $q$ 
in the right hand side of (\ref{ch1}).

The corresponding placements of $1, 2,\ldots,n$ and $q-n$ balls, 
henceforth called \textit{valid}, are the ones in which: 
\begin{enumit}
\item
$n$ occupies the first position from the left and
\item
at least $a$ balls separate an integer $k$ from the leftmost integer 
$i$ to the right of $k$, if $k > i$.
\end{enumit}
For example, the placement (\ref{pl}) is valid if $a \leq 2$. If a 
maximal string of consecutive balls has $p$ elements, we write 
$p = sa + r$ with $0 \leq r < a$ and think of the string as $s$ blocks 
consisting of $a$ balls each, simply refered to as \textit{a-blocks}, 
followed by $r$ balls. If $a=2$ then (\ref{pl}) has two $a$-blocks.

\smallskip
To construct the valid placements, let $j$ be the number of $a$-blocks. 
Place $j$ such blocks along a line and the integer $n$ first from the 
left, to guarantee ($i$). Insert $1, 2,\ldots,n-1$ in the $j$ spaces 
between the blocks and to the right of the last one, listing the 
integers within each space in increasing order to guarantee ($ii$). 
This can be done in $j^{n-1}$ ways. Finally, place the remaining 
$q - n -aj$ balls in the $n$ possible spaces between the integers and 
to the right of the last one, with at most $a-1$ in each space. The 
total number of ways is the coefficient of $y^{q-n}$ in (\ref{ch1}). 
This is clearly a polynomial in $q$ for $q > an$, hence (\ref{ch1}) 
holds specifically for all $q > an$.
\end{proof}

If $a=2$, the three-step procedure just described to construct (\ref{pl})
is the following:
\[ 4 \ \bigcirc \, \bigcirc \ \ \ \bigcirc \, \bigcirc \]
\[ 4 \ \bigcirc \, \bigcirc \ 2 \ 3 \ \bigcirc \, \bigcirc \ 1 \]
\[ 4 \ \bigcirc \, \bigcirc \, \bigcirc \ 2 \ 3 \ \bigcirc \, \bigcirc \ 
1 \ \bigcirc. \]
A simple application of Lemma 2.2 yields the more explicit formula 
for $\chi(\hat{\mathcal A}_n ^{[0,a]}, q)$, given in 
\cite[Thm.\ 9.7]{PS}.
\begin{corollary}
For all $a \geq 1$,
\[ \chi(\hat{\mathcal A}_n ^{[0,a]}, q) = \frac{q}{a^n} \
 S^n \, (1 + S + S^2 + \cdots + S^{a-1})^n \, q^{n-1}. \]
\end{corollary}
\begin{proof}
Formula (\ref{ch1}) can be written in the form
\[ \chi(\hat{\mathcal A}_n ^{[0,a]}, q) = \frac{q}{a^{n-1}} \
\sum_{k \equiv q-n \ ({\rm mod} \ a)} c_k \, (q-n-k)^{n-1}, \]
where the coefficients $c_k$ are as in Lemma 2.2. This lemma implies 
the proposed equality for $q > an$. Since both hand sides are polynomials 
in $q$, the equality follows for all $q$. 
\end{proof}

A similar formula follows for $\chi(\hat{\mathcal A}_n ^{[1,a]}, q)$. 
For convenience, as in \cite{Ath1,Ath2}, we use the notation 
$\tilde{\chi}({\mathcal A}, q) := \frac{1}{q} \ \chi({\mathcal A}, q)$.
\begin{proposition}{\rm (\cite[Thm.\ 4.3]{Ath1} \cite[Thm.\ 6.4.3]{Ath2})}
For all $a \geq 1$,
\[ \tilde{\chi}(\hat{\mathcal A}_n ^{[0,a]}, q) = 
   \tilde{\chi}(\hat{\mathcal A}_n ^{[1,a-1]}, q-n). \]
\end{proposition}
\begin{proof}
For $q$ large, the quantity on the right counts the $n$-tuples 
$(x_1, x_2,\ldots,x_n) \in {\mathbb Z}_{q-n} ^n$ which satisfy 
$x_i - x_j \neq 1,\ldots,a-1$ for all $1 \leq i < j \leq n$ and, 
say, $x_n = 0$. These $n$-tuples can be modeled again by placements 
of length $q-n$ of the integers $1, 2,\ldots,n$ and balls, in which 
more than one integer can occupy the same position since some of the 
$x_i$ may be equal. 

To define an explicit bijection with the valid placements of Proposition 
3.1, we start with a valid placement and remove a ball between any two 
consecutive integers, including the pair formed by the rightmost integer 
in the placement and $n$, which is the leftmost. If no ball lies between 
such a pair $(i, j)$ then we place $i$ and $j$ in the same position. For 
example, the placement (\ref{pl}) becomes 
\[ 4 \ \bigcirc \, \bigcirc \ \widehat{2 \ 3} \ \bigcirc \ 1 \] 
and corresponds to the $4$-tuple $(5, 3, 3, 0) \in {\mathbb Z}_6 ^4$. 
This map is clearly a bijection between the two kinds of placements. 
\end{proof}
\begin{corollary}{\rm (\cite[Thm.\ 9.7]{PS})}
For all $a \geq 1$,
\[ \chi(\hat{\mathcal A}_n ^{[1,a]}, q) = \frac{q}{(a+1)^n} \
 (1 + S + S^2 + \cdots + S^a)^n \, q^{n-1}. \] 
\qed
\end{corollary}

The special case $a=1$ of this corollary leads to another proof of 
Postnikov's theorem \cite{PS, St2}, initially conjectured by Linial
and Stanley. We give more details in Remark 1 of Section 6.

\section{Other root systems}

In this section we derive analogues of Corollaries 3.2 and 3.4 for the 
root systems $B_n, C_n, D_n$ and $BC_n$. The method we use follows 
closely that of Section 3. 

\vspace{0.1 in}
We need to adjust some of the terminology and reasoning of the previous 
section. Let $q$ be an odd positive integer. If 
$x = (x_1, x_2,\ldots,x_n)$ is an $n$-tuple of elements of 
${\mathbb Z}_q$ satisfying $x_i \neq 0$ for all $i$ and 
$x_i \neq \pm x_j$ for $i \neq j$, then we represent $x$ as a placement 
of the integers $1, 2,\ldots,n$, each with a $+$ or $-$ sign, and 
$\frac{q-1}{2} - n$ indistinguishable balls along a line, with an 
extra zero in the first position from the left. For example, omitting 
the $+$ signs, for $q = 27$ and $n = 6$ we have the placement
\begin{equation}
0 \ \bigcirc \ 2 \ \bigcirc \ 3 \ -5 \ \bigcirc \, \bigcirc \ 4 \ -1 \ 
\bigcirc \, \bigcirc \ -6 \ \bigcirc.
\label{pl2}
\end{equation}
Such a placement corresponds to the $n$-tuple $x$ for which $x_i + 1$ 
or $-x_i + 1$ is the position that $i$ or $-i$ occupies, respectively, 
counting from the left. The placement (\ref{pl2}) corresponds to the 
$6$-tuple $(-9, 2, 4, 8, -5, -12)$ of elements of ${\mathbb Z}_{27}$. 

\vspace{0.1 in}
We first derive the analogues of Corollary 3.2 in the four cases of 
interest. The symbol $\prec$ refers to the total order of the integers 
\[ 1 \prec 2 \prec 3 \prec \cdots \prec 0 \prec \cdots \prec -3 \prec 
-2 \prec -1. \]

\vspace{0.1 in}
\textit{The root system $BC_n$}. Recall that $\hat{\mathcal {BC}}_n 
^{[0,a]}$ has hyperplanes
\begin{equation}
\begin{tabular}{l}
$x_i = 0, 1,\ldots,a \ \ \for \ \ 1 \leq i \leq n$,\\
$2x_i = 0, 1,\ldots,a \ \ \for \ \ 1 \leq i \leq n$,\\
$x_i - x_j = 0, 1,\ldots,a \ \ \for \ \ 1 \leq i < j \leq n$,\\
$x_i + x_j = 0, 1,\ldots,a \ \ \for \ \ 1 \leq i < j \leq n$.
\end{tabular}
\label{arrBC1}
\end{equation}
\begin{proposition}
For $a \geq 1$, $\chi(\hat{\mathcal {BC}}_n ^{[0,a]}, q)$ is equal to
\[ \frac{2}{a^{n+1}} \ S^{2n+1} \, (1 + S^2 + S^4 + \cdots + S^{2a-2})^n \, 
  (1 + S^2 + S^4 + \cdots + S^{a-2}) \ q^n \]
if $a$ is even and
\[ \frac{1}{a^{n+1}} \ S^{2n+1} \, (1 + S^2 + S^4 + \cdots + S^{2a-2})^n \, 
  (1 + S + S^2 + \cdots + S^{a-1}) \ q^n \]
if $a$ is odd.
\end{proposition}
\begin{proof}
By Theorem 2.1, for sufficiently large odd $q$, 
$\chi(\hat{\mathcal {BC}}_n ^{[0,a]}, q)$ counts the number of 
$n$-tuples $x = (x_1, x_2,\ldots,x_n) \in {\mathbb Z}_q ^n$ for which 
none of the equalities (\ref{arrBC1}) holds in ${\mathbb Z}_q$. The 
corresponding placements of integers and balls are the ones in which:
\begin{enumit}
\item
at least $a$ balls are placed between $0$ and the leftmost nonzero 
integer, if this integer is positive,
\item
at least $\lfloor \frac{a+1}{2} \rfloor$ balls are placed to the right 
of the rightmost integer, if this integer is negative and
\item
at least $a$ balls separate a nonzero integer $k$ from the leftmost 
integer $i$ to the right of $k$ if $k \succ i$.
\end{enumit}
We call again these placements \textit{valid}. The placement (\ref{pl2}) 
is valid for $a=1$ but it is not for $a \geq 2$. Conditions $(i)$ and 
$(ii)$ guarantee that no equation of the first two kinds in 
(\ref{arrBC1}) holds. For example $2x_i \neq 1$, or equivalently 
$-x_i \neq \frac{q-1}{2}$, requires that the last position from the 
right is not occupied by $-i$. Condition $(iii)$ takes care of the 
remaining two kinds of equations. 

\smallskip
To construct the valid placements, place $j$ $a$-blocks along a line, 
as in the proof of Proposition 3.1, and $0$ to the left. Insert 
$1, 2,\ldots,n$, each with a sign, in one of the $j+1$ possible spaces 
between $0$ and the $a$-blocks and to the right of the last $a$-block. 
List the integers within each space in increasing order with respect 
to $\prec$, to guarantee ($iii$), and force the $-$ sign in the space 
immediately to the right of $0$, to guarantee ($i$). Then distribute 
the remaining $\frac{q-1}{2} - n -aj$ balls between the integers, in 
blocks of at most $a-1$. To take care of ($ii$), we distinguish two 
cases according to whether there is a negative integer to the right
of the rightmost $a$-block or not. It follows that
\[ \chi(\hat{\mathcal {BC}}_n ^{[0,a]}, q) = [y^{p-n}] \ \,
 (\phi_a (y))^{n+1} \, \sum_{j=0}^{\infty} \, (2j)^n \, y^{aj} \ \ + \]
\vspace*{-0.18 in}
\[ [y^{p-n}] \ \left( y^{\lfloor \frac{a+1}{2} \rfloor} + \cdots + 
 y^{a-1} \right) \, (\phi_a (y))^n \, \sum_{j=0}^{\infty} \, 
((2j+1)^n - (2j)^n) \, y^{aj}, \]
where $p = \frac{q-1}{2}$. The quantity $(2j+1)^n - (2j)^n$ in the 
second summand stands for the number of ways to insert the integers 
$1, 2,\ldots,n$ with signs in $j+1$ possible spaces with the $-$ sign 
forced in the first space and at least one $-$ sign in the last. 

\smallskip
We now extract the coefficients of $y^{p-n}$ and use Lemma 2.2 as in 
Corollary 3.2 to get the proposed expressions, after some 
straightforward algebraic manipulations. Note that $(2j+1)^n - (2j)^n$ 
has degree $n-1$ in $j$, so Lemma 2.2 applies to the second summand 
as well.
\end{proof}

\vspace{0.1 in}
The derivations in the other three cases involve some complications 
but are treated in a similar way, so we will omit most of the details.
We let $p = \frac{q-1}{2}$ until the end of this section.  

\vspace{0.1 in}
\textit{The root system $C_n$}. The arrangement $\hat{\mathcal C}_n 
^{[0,a]}$ lacks the first set of hyperplanes in (\ref{arrBC1}). 
\begin{proposition}
For $a \geq 1$, $\chi(\hat{\mathcal C}_n ^{[0,a]}, q)$ is equal to
\[ \frac{4}{a^{n+1}} \ S^{2n+1} \, (1 + S^2 + S^4 + \cdots + 
 S^{2a-2})^{n-1} \, (1 + S^2 + S^4 + \cdots + S^{a-2})^2 \ q^n \]
if $a$ is even and
\[ \frac{1}{a^{n+1}} \ S^{2n} \, (1 + S^2 + S^4 + \cdots + 
S^{2a-2})^{n-1} \, (1 + S + S^2 + \cdots + S^{a-1})^2 \ q^n \]
if $a$ is odd.
\end{proposition}
\begin{proof}
The valid placements in this case are as for $BC_n$ except that, in 
condition $(i)$, $a$ is replaced by $\lfloor \frac{a}{2} \rfloor$. To 
count these placements we now distinguish four cases, according to 
whether there is a positive integer between zero and the leftmost 
$a$-block and whether there is a negative integer to the right 
of the rightmost $a$-block. It follows that
\[ \chi(\hat{\mathcal C}_n ^{[0,a]}, q) = [y^{p-n}] \
 (\phi_a (y))^{n+1} \, \sum_{j=0}^{\infty} \, (2j)^n \, y^{aj} \ \ + \]
\vspace*{-0.18 in}
\[ [y^{p-n}] \ \left( y^{\lfloor \frac{a}{2} \rfloor} + \cdots + y^{a-1} 
 \right) \, (\phi_a (y))^n \, \sum_{j=0}^{\infty} \, ((2j+1)^n - (2j)^n) 
 \, y^{aj} \ \ + \]
\vspace*{-0.08 in}
\[ [y^{p-n}] \ \left( y^{\lfloor \frac{a+1}{2} \rfloor} + \cdots + y^{a-1} 
 \right) \, (\phi_a (y))^n \, \sum_{j=0}^{\infty} \, ((2j+1)^n - (2j)^n) \, 
y^{aj} \ \ + \]
\vspace*{-0.18 in}
\[ [y^{p-n}] \ \left( y^{\lfloor \frac{a}{2} \rfloor} + \cdots + y^{a-1} 
 \right) \left( y^{\lfloor \frac{a+1}{2} \rfloor} + \cdots + y^{a-1} \right) 
 \, (\phi_a (y))^{n-1} \, \sum_{j=0}^{\infty} \, a_{2j+2} \, y^{aj}, \]
where
\begin{equation}
a_j = j^n - 2(j-1)^n + (j-2)^n.
\label{a_j}
\end{equation}
The result follows in a straightforward way, as before. Note that the 
degree of $a_{2j+2}$ in $j$ is at most $n-2$ and hence Lemma 2.2 
applies to the last summand as well.
\end{proof}

\vspace{0.1 in}
\textit{The root system $B_n$}. The arrangement $\hat{\mathcal B}_n 
^{[0,a]}$ lacks the second set of hyperplanes in (\ref{arrBC1}). The 
proof of the following proposition is indirect. 
\begin{proposition}
For $a \geq 1$,
\[ \chi(\hat{\mathcal B}_n ^{[0,a]}, q) = 
   \chi(\hat{\mathcal C}_n ^{[0,a]}, q). \]
\end{proposition}
\begin{proof}
Let $l, m$ denote the last two integers in a placement and $s, t$ the 
number of balls between $l$ and $m$ and to the right of $m$, 
respectively. For the placement (\ref{pl2}) we have $l = -1$, $m = -6$, 
$s = 2$ and $t = 1$. The valid placements for $\hat{\mathcal B}_n 
^{[0,a]}$ are the ones which satisfy conditions $(i)$ and $(iii)$ of 
the $BC_n$ case and also:

\smallskip
\hspace{0.15 in} $(ii^{\prime})$ $2s + t \geq a-1$ if $l \succ -m$. 

\smallskip
\noindent
Indeed, the conditions $x_i \pm x_j \neq 0, 1,\ldots,a$ require that 
$(iii)$ holds, with the extra assumption $k \neq -i$, 
if we extend a placement, say (\ref{pl2}), to the rest of the classes 
mod $q$ as  
\[ 0 \ \bigcirc \ 2 \ \bigcirc \ 3 \ -5 \ \bigcirc \, \bigcirc \ 4 
\ -1 \ \bigcirc \, \bigcirc \ -6 \ \bigcirc \ \bigcirc \ 6 \ \bigcirc 
\, \bigcirc \ 1 \ -4 \ \cdots. \]
This also implies $(ii^{\prime})$. Note that $(ii^{\prime})$ 
follows from $(iii)$ if $l \succ m$ but is essential otherwise. It 
is redundant in the cases of $BC_n$ and $C_n$ because of $(ii)$. To 
count the valid placements in this case, we first count those which 
satisfy $(i)$ and $(iii)$ and then subtract the ones which violate 
$(ii^{\prime})$. For large odd $q$, it follows that 
$\chi(\hat{\mathcal B}_n ^{[0,a]}, q)$ is the coefficient of $y^{p-n}$ 
in the expression
\[ (\phi_a (y))^{n+1} \, \sum_{j=0}^{\infty} \, (2j+1)^n \, y^{aj} - 
 f_{a-2}(y) \ (\phi_a (y))^{n-1} \, \sum_{j=0}^{\infty} \, a_j^{\prime} 
 \, y^{aj}, \]
where
\[ f_k(y) := \sum_{\scriptsize \begin{array}{c} s, t \geq 0 \\ 2s+t 
\leq k \end{array}} y^{s+t} \]
and $a_j^{\prime}$ is the number of ways to insert 
$1, 2,\ldots,n$ with signs in $j+1$ spaces and list the 
integers in each space in increasing order with respect to $\prec$ so 
that the last two integers $l, m$ appear in the last space and satisfy 
$l \succ -m$, in addition to $l \prec m$. It is easy to check that  
\begin{equation}
f_k(y) = 1 + 2y + 3y^2 + \cdots + 2y^{k-1} + y^k = 
\left\{ \begin{array}{ll} (\phi_{r+1}(y))^2, & \mbox{$k=2r$};\\
         \phi_r(y) \, \phi_{r+1}(y), & \mbox{$k=2r-1$} \end{array} 
\right.
\label{f}
\end{equation}
and that $a_j^{\prime} = \sum_{k=2}^{n} {n \choose k} (2^k - 2) \,
(2j-1)^{n-k} = a_{2j+1}$, defined by (\ref{a_j}). In this sum, $k$
stands for the number of integers in the last space. Using Lemma 2.2 
as before, we arrive at the same expression for 
$\chi(\hat{\mathcal B}_n ^{[0,a]}, q)$ as the one obtained earlier 
for $\chi(\hat{\mathcal C}_n ^{[0,a]}, q)$.
\end{proof}

\vspace{0.1 in}
\textit{The root system $D_n$}. The arrangement $\hat{\mathcal D}_n 
^{[0,a]}$ lacks the first two sets of hyperplanes in (\ref{arrBC1}). 
Let ${\mathcal Q}_n$ be the arrangement of coordinate hyperplanes 
$x_i = 0$ in ${\mathbb R}^n$. Then $\hat{\mathcal D}_n ^{[0,a]} \cup 
{\mathcal Q}_n$ has hyperplanes
\begin{equation}
\begin{tabular}{l}
$x_i = 0 \ \ \for \ \ 1 \leq i \leq n$,\\
$x_i - x_j = 0, 1,\ldots,a \ \ \for \ \ 1 \leq i < j \leq n$,\\
$x_i + x_j = 0, 1,\ldots,a \ \ \for \ \ 1 \leq i < j \leq n$.
\end{tabular}
\label{arrDQ}
\end{equation}
We first prove the following lemma.
\begin{lemma}
For $a \geq 1$ and $n \geq 3$, $\chi(\hat{\mathcal D}_n ^{[0,a]} \cup 
{\mathcal Q}_n, q)$ is equal to
\[ \frac{4 S^{2n-1}}{a^{n+1}} \ \left(\phi_a(S^2)\right)^{n-3}
\, \left(\phi_{a/2}(S^2)\right)^4 \, (1 + 3S^2 - S^a + S^{a+2}) \ q^n \]
if $a$ is even and
\[ \frac{S^{2n-1}}{a^{n+1}} \ \left(\phi_a(S^2)\right)^{n-3}
  \, (\phi_a(S))^4 \ \frac{2-S^{a-1}+S^a}{1+S} \ q^n \]
if $a$ is odd.
\end{lemma}
\begin{proof}
Let $l^{\prime}, m^{\prime}$ denote the first two integers in a 
placement and $s^{\prime}, t^{\prime}$ the number of balls to the 
left of $l^{\prime}$ and between $l^{\prime}$ and $m^{\prime}$, 
respectively. For the placement (\ref{pl2}) we have $l^{\prime} = 2$, 
$m^{\prime} = 3$ and $s^{\prime} = t^{\prime}= 1$. The valid placements 
for $\hat{\mathcal D}_n ^{[0,a]} \cup {\mathcal Q}_n$ are the ones which 
satisfy conditions $(ii^{\prime})$ and $(iii)$ of the $B_n$ case (see 
the proof of Proposition 4.1 for $(iii)$) and also:

\smallskip
\hspace{0.15 in} $(i \, ^{\prime})$ $2s^{\prime} + t^{\prime} \geq a-2$ 
if $-l^{\prime} \succ m^{\prime}$.

\smallskip
\noindent
This is implied by $(iii)$ if $l^{\prime} \succ m^{\prime}$ but is 
essential otherwise. It is redundant in the cases of $BC_n$, $C_n$ 
and $B_n$ because of $(i)$ and its $C_n$ analogue. We count these valid 
placements as in the $B_n$ case, using a simple inclusion-exclusion to 
handle both $(i^{\prime})$ and $(ii^{\prime})$. It follows that, for 
large odd $q$, $\chi(\hat{\mathcal D}_n ^{[0,a]} \cup {\mathcal Q}_n, q)$ 
is the coefficient of $y^{p-n}$ in the expression
\[ (\phi_a (y))^{n+1} \, \sum_{j=0}^{\infty} \, (2j+2)^n \, y^{aj} \ - \ 
f_{a-3}(y) \ (\phi_a (y))^{n-1} \, \sum_{j=0}^{\infty} \, a_{2j+2} \, 
y^{aj} \ - \]
\vspace*{-0.12 in}
\[ - \ f_{a-2}(y) \ (\phi_a (y))^{n-1} \, \sum_{j=0}^{\infty} \, a_{2j+2} 
\, y^{aj} \ \ + \ \ f_{a-3}(y) \, f_{a-2}(y) \ (\phi_a (y))^{n-3} \,
 \sum_{j=0}^{\infty} \, b_j \, y^{aj}, \]
where we have used the notation in (\ref{a_j}) and (\ref{f}) and 
\[ b_j = (2j+2)^n - 4(2j+1)^n + 6(2j)^n - 4(2j-1)^n + (2j-2)^n, \]
by a computation similar to the one for $a_j^{\prime}$ in the proof
of Proposition 4.3. We extract this coefficient and factor the 
resulting expression appropriately to get the result.
\end{proof}

We now compute $\chi(\hat{\mathcal D}_n ^{[0,a]}, q)$ for $n \geq 3$. 
It is easy to check that $\chi(\hat{\mathcal D}_2 ^{[0,a]}, q) = 
(q-a-1)^2$ for all $a$.
\begin{proposition}
For $a \geq 1$ and $n \geq 3$, $\chi(\hat{\mathcal D}_n ^{[0,a]}, q)$ 
is equal to
\[ \frac{8 S^{2n-1}}{a^{n+1}} \ (1 + S^2) \, 
 (1 + S^2 + S^4 + \cdots + S^{2a-2})^{n-3} \, 
 (1 + S^2 + S^4 + \cdots + S^{a-2})^4 \ q^n \]
if $a$ is even and
\[ \frac{1}{a^{n+1}} \ S^{2n-2} \, (1 + S^2 + S^4 + \cdots + 
S^{2a-2})^{n-3} \, (1 + S + S^2 + \cdots + S^{a-1})^4 \ q^n \]
if $a$ is odd.
\end{proposition}
\begin{proof}
By Theorem 2.1, for large odd $q$, $\chi(\hat{\mathcal D}_n ^{[0,a]}, q)$ 
counts the number of $n$-tuples $x = (x_1, x_2,\ldots,x_n) \in 
{\mathbb Z}_q ^n$ which satisfy 
\begin{equation}
x_i \pm x_j \neq 0, 1,\ldots,a
\label{pm}
\end{equation} 
in ${\mathbb Z}_q$ for all $1 \leq i < j \leq n$. The ones which also 
satisfy $x_i \neq 0$ for all $i$ were counted in the previous lemma. 
Therefore, the characteristic polynomial of $\hat{\mathcal D}_n ^{[0,a]}$ 
is the sum of that of $\hat{\mathcal D}_n ^{[0,a]} \cup {\mathcal Q}_n$ 
and $\psi(q)$, where $\psi(q)$ is the number of $n$-tuples $x$ for which 
(\ref{pm}) holds and $x_i = 0$ for at least one, and hence exactly one 
$i$. These can be modeled by placements which satisfy conditions 
$(ii^{\prime})$ and $(iii)$ of the $B_n$ case but have a negative 
integer in the leftmost position, instead of $0$. For example, 
\[ -2 \ \bigcirc \ 3 \ -5 \ \bigcirc \, \bigcirc \ 4 \ -1 \ 
\bigcirc \, \bigcirc \ -6 \ \bigcirc \]
corresponds to the $6$-tuple $(-7, 0, 2, 6, -3, -10) \in 
{\mathbb Z}_{23} ^6$. Thus, when constructing these placements, 
at least one negative integer is inserted to the left of the
leftmost $a$-block but no positive one. The argument in the $B_n$ 
case shows that 
\[ \psi(q) = [y^{p-n+1}] \ 
 (\phi_a (y))^n \, \sum_{j=0}^{\infty} \, ((2j+1)^n - (2j)^n) \, 
y^{aj}\ \ - \]
\vspace*{-0.18 in}
\[ [y^{p-n+1}] \ f_{a-2}(y) \ (\phi_a (y))^{n-2} \,
 \sum_{j=0}^{\infty} \, d_j \, y^{aj}, \]
where
\[ d_j = a_{2j+1} - a_{2j} = (2j+1)^n - 3(2j)^n + 3(2j-1)^n - 
(2j-2)^n. \]
It follows that $\psi(q)$ is equal to
\[ \frac{4S^{2n-1}}{a^{n+1}} \ (1 - S^a) \, (\phi_a (S^2))^{n-2} \, 
 (\phi_{a/2} (S^2))^2 \ q^n \]
if $a$ is even and
\[ \frac{S^{2n-2}}{a^{n+1}} \  (1 - S^a) \, (\phi_a (S^2))^{n-2} \,
 (\phi_{a/2} (S^2))^2 \ q^n \]
if $a$ is odd. These expressions and Lemma 4.4 imply the result.
\end{proof}

\vspace{0.1 in}
The analogue of Proposition 3.3 was derived for most of the cases of 
interest in \cite{Ath2}. 
\begin{proposition}{\rm (\cite[Thm.\ 7.2.4 and Thm.\ 7.2.7]{Ath2})}
If $\Phi = B_n$ or $D_n$ and $a \geq 1$ or $\Phi = C_n$ or $BC_n$ and 
$a \geq 2$ is even, then
\[ \chi\left(\hat{\mathcal A} ^{[0,a]} \, (\Phi), q\right) = 
   \chi\left(\hat{\mathcal A} ^{[1,a-1]} \, (\Phi), q-h\right), \]
where
\[ h = \left\{ 
\begin{tabular}{l} 
$2n-2$, \textrm{if $\Phi = D_n$;}\\
$2n$, \textrm{otherwise.} 
\end{tabular} \right. \]
\end{proposition}
\begin{proof}
For large odd $q$, the quantities on the right hand side count the 
$n$-tuples $(x_1, x_2,\ldots,x_n) \in {\mathbb Z}_{q-h} ^n$ which 
satisfy $x_i \pm x_j \neq 1,\ldots,a-1$ for all $1 \leq i < j \leq n$ 
and some of the conditions $x_i \neq 1,\ldots,a-1$ and 
$2x_i \neq 1,\ldots,a-1$, depending on the case. These $n$-tuples can 
be modeled by placements of length $\frac{q+1}{2}$, as described in 
the beginning of this section, except that more than one integer can 
occupy the same position, possibly the leftmost, labeled with a zero 
otherwise.

\smallskip
In each case there is an explicit bijection with the valid placements 
of Propositions 4.1 -- 4.5. Given a valid placement, we remove a ball 
between any two consecutive integers, as in the proof of Proposition 3.3. 
These pairs of integers include the one formed by $0$ and the leftmost 
nonzero integer in the cases of $B_n, C_n$ and $BC_n$ but not in the 
case of $D_n$. Also, in all four cases we leave the number of balls to 
the right of the rightmost integer unchanged. For example, the placement 
(\ref{pl2}) becomes 
\[ 0 \ \bigcirc \ 2 \ \ \widehat{3 -\hspace{-1mm}5} \ \bigcirc \ 
\widehat{4 -\hspace{-1mm}1} \ \bigcirc \ -6 \ \ \bigcirc \]
in the case of $D_n$ and
\[ 0 \ \ 2 \ \ \widehat{3 -\hspace{-1mm}5} \ \bigcirc \ \widehat{4 
 -\hspace{-1mm}1} \ \bigcirc \ -6 \ \ \bigcirc \]
in the three other cases. They correspond to the $6$-tuples 
$(-5, 2, 3, 5, -3, -7) \in {\mathbb Z}_{17} ^6$ and 
$(-4, 1, 2, 4, -2, -6) \in {\mathbb Z}_{15} ^6$ respectively. It is 
easy to see that this map is indeed a bijection in each case.
\end{proof}

The bijection just described breaks down in the cases of odd $a$ for 
$\Phi = C_n$ or $BC_n$, which need special care. The following 
proposition was conjectured in \cite{Ath2}.
\begin{proposition}{\rm (\cite[Conjecture 7.2.8]{Ath2})}
For all odd $a \geq 1$, 
\[ \chi(\hat{\mathcal C}_n ^{[0,a]}, q) = 
   \chi(\hat{\mathcal C}_n ^{[1,a-1]}, q-2n) \]
and 
\[ \chi(\hat{\mathcal {BC}}_n ^{[0,a]}, q) = 
   \chi(\hat{\mathcal {BC}}_n ^{[1,a-1]}, q-2n-1). \]
\end{proposition}
\begin{proof}
For the first statement, let $q$ be a large odd integer. Start with a 
valid placement, as described in the proof of Proposition 4.2. Read it 
from left to right, switch the $+$ signs to $-$ and vice versa and 
disregard $0$, to get a new placement. Finally remove a ball between 
consecutive integers, as in the proof of Proposition 4.6, but leave 
the number of balls in the far left and far right unchanged, to get 
a placement counted by the right hand side. For example, (\ref{pl2}) 
becomes
\[ \bigcirc \ 6 \ \bigcirc \ \widehat{1 -\hspace{-1mm}4} \ \bigcirc \ 
\widehat{5 -\hspace{-1mm}3} \ -2 \ \bigcirc, \]
which corresponds to the $6$-tuple $(3, -6, -5, -3, 5, 1)$ of elements 
of ${\mathbb Z}_{15} ^6$. It is easy to check that this map is a 
bijection.

\smallskip
Note that a direct bijective proof by Theorem 2.1 is not possible for 
the second statement since $q$ and $q-2n-1$ cannot both be odd. Once 
the valid placements for $\hat{\mathcal {BC}}_n ^{[1,a-1]}$ are 
described explicitly, an argument similar to the one in the proof of 
Proposition 4.1 shows that
\[ \chi(\hat{\mathcal {BC}}_n ^{[1,a-1]}, q) = [y^p] \ \,
 (\phi_a (y))^{n+1} \, \sum_{j=0}^{\infty} \, (2j)^n \, y^{aj} \ \ + \]
\vspace*{-0.18 in}
\[ [y^p] \ \left( y^{\frac{a-1}{2}} + \cdots + y^{a-1} \right) \, 
   (\phi_a (y))^n \, \sum_{j=0}^{\infty} \, ((2j+1)^n - (2j)^n) \, 
y^{aj}. \]
This implies the result indirectly, by 
comparison to the formula of Proposition 4.1.
\end{proof}

Analogues of Corollary 3.4 follow in all four cases. For example, in 
the case of $BC_n$ we have the following corollary.   
\begin{corollary}
For all $a \geq 1$, $\chi(\hat{\mathcal {BC}}_n ^{[1,a]}, q)$ is equal to
\[ \frac{2 S}{(a+1)^{n+1}} \ (1 + S^2 + S^4 + \cdots + S^{2a})^n \, 
  (1 + S^2 + S^4 + \cdots + S^{a-1}) \ q^n \]
if $a$ is odd and
\[ \frac{1}{(a+1)^{n+1}} \ (1 + S^2 + S^4 + \cdots + S^{2a})^n \, 
  (1 + S + S^2 + \cdots + S^a) \ q^n \]
if $a$ is even.
\qed
\end{corollary}

\section{Proof of the Main Theorem}

The results of Sections 3 and 4 imply a crucial case of Theorem 1.2 
via the following lemma. This lemma was used by Postnikov and Stanley 
in \cite{PS} to prove Conjecture 1.1 for the root system $A_{n-1}$.
\begin{lemma}{\rm (\cite[Lemma 9.12]{PS})}
If $g, f \in {\mathbb C}[q]$ are such that $g$ has degree $d$, all 
roots of $g$ have absolute value $1$ and all roots of $f$ have real 
part equal to $r$, then all roots of $g(S) f$ have real part equal 
to $r + d/2$.
\end{lemma}
\begin{corollary}
Conjecture 1.1 holds for ${\mathcal A} = \hat{\mathcal A}^{[0, b]} \, 
(\Phi)$, $\hat{\mathcal A}^{[1, b]} \, (\Phi)$ if $b$ is a positive 
integer and $\Phi$ is one of $A_{n-1}$, $B_n$, $C_n$, $D_n$ or $BC_n$ 
for some $n \geq 2$.
\end{corollary}
\begin{proof}
Combine the results of Sections 3 and 4 with Lemma 5.1.
\end{proof}

To complete the proof of Theorem 1.2 we need one last result. The 
first statement in the following proposition is the content of 
\cite[Thm.\ 7.2.1]{Ath2}. We note that the argument in the case of 
$C_n$, given there, was oversimplified.
\begin{proposition}
Let $a, b$ be integers satisfying $0 \leq a \leq b$. If $\Phi$ is one 
of $A_{n-1}$, $B_n$, $C_n$ or $D_n$ then
\[ \chi(\hat{\mathcal A}^{[-a,b]} \, (\Phi), q) = 
   \chi(\hat{\mathcal A}^{[0,b-a]} \, (\Phi), q-ah), \]
where
\[ h = \left\{ 
\begin{tabular}{l} 
$n$, \textrm{if $\Phi = A_{n-1}$;}\\
$2n$, \textrm{if $\Phi = B_n$ or $C_n$;}\\ 
$2n-2$, \textrm{if $\Phi = D_n$.}
\end{tabular} \right. \]
For $\Phi = BC_n$,
\[ \chi(\hat{\mathcal {BC}}_n ^{[-a,b]}, q) = 
   \left\{ 
\begin{tabular}{l} 
$\chi\left(\hat{\mathcal {BC}}_n ^{[0,b-a]}, q-(2n+1)a-1\right)$, 
 \textrm{if both $a$ and $b$ are odd;}
 \\
$\chi\left(\hat{\mathcal {BC}}_n ^{[0,b-a]}, q-(2n+1)a\right)$, 
 \textrm{otherwise.}
\end{tabular} \right. \]
\end{proposition}
\begin{proof}
Let $\Phi$ be as above but assume that either $a$ is even or $b$ is 
odd if $\Phi = C_n$ or $BC_n$. For large $q$ if $\Phi = A_{n-1}$ and 
large odd $q$ otherwise, both hand sides of the proposed equalities 
count placements of a certain kind. To obtain a bijection, we start 
with a placement counted by the right hand side and simply add $a$ 
balls between consecutive integers, as defined in Propositions 3.3 
and 4.6, except that we only add $\lfloor \frac{a}{2} \rfloor$ balls 
immediately to the left of $0$ if $\Phi = C_n$ and that we add 
$\lfloor \frac{a+1}{2} \rfloor$ balls to the right of the rightmost 
integer if $\Phi = C_n$ or $BC_n$. If $a=1$ then (\ref{pl}) becomes
\[ 4 \ \bigcirc \, \bigcirc \, \bigcirc \, \bigcirc \ 2 \ \bigcirc \ 
3 \ \bigcirc  \, \bigcirc \, \bigcirc \ 1 \ \bigcirc \, \bigcirc \]
and (\ref{pl2}) becomes
\[ 0 \ \bigcirc \, \bigcirc \ 2 \ \bigcirc \, \bigcirc \ 3 \ \bigcirc 
\ -5 \ \bigcirc \, \bigcirc \, \bigcirc \ 4 \ \bigcirc \ -1 \ \bigcirc 
\, \bigcirc \, \bigcirc \ -6 \ \, \bigcirc, \]
\[ 0 \ \bigcirc \ 2 \ \bigcirc \, \bigcirc \ 3 \ \bigcirc \ -5 \ 
\bigcirc \, \bigcirc \, \bigcirc \ 4 \ \bigcirc \ -1 \ \bigcirc \, 
\bigcirc \, \bigcirc \ -6 \ \, \bigcirc, \]
\[ 0 \ \bigcirc \ 2 \ \bigcirc \, \bigcirc \ 3 \ \bigcirc \ -5 \ 
\bigcirc \, \bigcirc \, \bigcirc \ 4 \ \bigcirc \ -1 \ \bigcirc \, 
\bigcirc \, \bigcirc \ -6 \ \bigcirc \, \bigcirc \]
or
\[ 0 \ \bigcirc \, \bigcirc \ 2 \ \bigcirc \, \bigcirc \ 3 \ \bigcirc 
\ -5 \ \bigcirc \, \bigcirc \, \bigcirc \ 4 \ \bigcirc \ -1 \ \bigcirc 
\, \bigcirc \, \bigcirc \ -6 \ \bigcirc \, \bigcirc, \]
if $\Phi = B_n$, $D_n$, $C_n$ or $BC_n$, respectively. This map is 
easily seen to be a bijection in each case.

\smallskip
Now suppose that $a$ is odd, $b$ is even and $\Phi = C_n$ or $BC_n$. 
The map described above fails to be well defined in these cases. 
Moreover, for $\Phi = BC_n$, a direct bijective proof is not possible 
since $q$ and $q-(2n+1)a$ cannot both be odd. One way to overcome this 
difficulty is to prove instead that
\[ \chi(\hat{\mathcal {BC}}_n ^{[-a,b]}, q) = 
   \chi\left(\hat{\mathcal {BC}}_n ^{[1,b-a-1]}, q-(2n+1)(a+1)\right) \]
and
\[ \chi(\hat{\mathcal C}_n ^{[-a,b]}, q) = 
   \chi\left(\hat{\mathcal C}_n ^{[1,b-a-1]}, q-2n(a+1)\right), \]
which are equivalent to the desired formulae by Proposition 4.7. Note 
that we have the empty arrangement in ${\mathbb R}^n$ on the right if 
$b=a+1$. Now a bijective proof is possible. Start with a placement 
counted by the left hand side and remove $a+1$ balls between consecutive 
integers, but only $\frac{a+1}{2}$ in the far right if $\Phi = BC_n$ 
and in the far left and far right if $\Phi = C_n$.
\end{proof}

\vspace{0.1 in}
We now return to the proof of the main thorem.

\vspace{0.1 in}
\noindent
\textit{Proof of Theorem 1.2.} Combine Corollary 5.2 and 
Proposition 5.3.
\qed

\section{Remarks}

1. For $a=1$, Corollary 3.4 yields the expression 
\[ \frac{q}{2^n} \sum_{k=0}^{n} {n \choose k} \, (q - k)^{n-1} \]
for the characteristic polynomial of the Linial arrangement of
hyperplanes in ${\mathbb R}^n$ $x_i - x_j = 1$ for $i < j$. It follows
via Zaslavsky's theorem \cite{Za} that the number $g_n$ of regions 
into which this arrangement dissects ${\mathbb R}^n$ is
\[ \frac{1}{2^n} \sum_{k=0}^{n} {n \choose k} \, (k + 1)^{n-1}, \]
which is also the number $f_n$ of \textit{alternating trees} on $n+1$ 
vertices \cite{Po}. The fact that $g_n = f_n$ was conjectured by 
Linial and Stanley and first proved by Postnikov \cite{PS, St2}. No 
bijective proof of this fact is known.

\vspace{0.1 in}
2. The results of Section 4 yield similar expressions for the number 
of regions of the Linial arrangement $\hat{\mathcal A} ^{[1]} \, 
(\Phi)$ for $\Phi = B_n, C_n, D_n$ and $BC_n$. This expression is
\[ 2 \sum_{k=0}^{n-1} {n-1 \choose k} \, (k + 1)^n \]
if $\Phi = B_n$ or $C_n$,
\[ 4 \sum_{k=0}^{n-2} {n-2 \choose k} \, (k + 1)^n \]
if $\Phi = D_n$ and
\[ \sum_{k=0}^{n} {n \choose k} \, (k + 1)^n \]
if $\Phi = BC_n$. It would be interesting to find combinatorial 
interpretations to these numbers similar to the one in the case
of $A_{n-1}$.

\vspace{0.1 in}
3. The reasoning in Section 4 can be applied to the more general family 
of deformations of the form
\begin{equation}
\begin{tabular}{l}
$x_i = 0, 1, 2,\ldots,b \ \ \for \ \ 1 \leq i \leq n$,\\
$2x_i = 1, 3,\ldots,2c-1 \ \ \for \ \ 1 \leq i \leq n$,\\
$x_i - x_j = 0, 1,\ldots,a \ \ \for \ \ 1 \leq i < j \leq n$,\\
$x_i + x_j = 0, 1,\ldots,a \ \ \for \ \ 1 \leq i < j \leq n.$
\end{tabular}
\label{arr1}
\end{equation}
We will only mention the special case $a = b = c$. The resulting 
arrangement is not one of the deformations of interest but the 
following proposition implies via Lemma 5.1 that Conjecture 1.1 still 
holds in this case and suggests that the conjecture is true in an even more 
general setting. Furthermore, the corresponding formula is easier to 
obtain.
\begin{proposition}
For $a = b= c \geq 1$, the arrangement (\ref{arr1}) has characteristic 
polynomial
\[ \frac{1}{a^{n+1}} \ 
   S^{2n+1} \, (1 + S^2 + S^4 + \cdots + S^{2a-2})^{n+1} \, q^n. \]
Also, the arrangement
\begin{equation}
\begin{tabular}{l}
$2x_i = 1, 2,\ldots,2a-1 \ \ \for \ \ 1 \leq i \leq n$,\\
$x_i - x_j = 1, 2,\ldots,a-1 \ \ \for \ \ 1 \leq i < j \leq n$,\\
$x_i + x_j = 1, 2,\ldots,a-1 \ \ \for \ \ 1 \leq i < j \leq n$
\end{tabular}
\label{arr2}
\end{equation}
has characteristic polynomial
\[ \frac{1}{(a+1)^{n+1}} \
   S \, (1 + S^2 + S^4 + \cdots + S^{2a})^{n+1} \, q^n. \]
\end{proposition}
\begin{proof}
The argument in the proof of Proposition 4.1 yields the expression
\[ [y^{p-n}] \ \, (1 + y + y^2 + \cdots + y^{a-1})^{n+1} \, 
 \sum_{j=0}^{\infty} \, (2j)^n \, y^{aj} \]
for the characteristic polynomial of the first arrangement in question. 
This implies the proposed formula, as well as the formula for the 
second arrangement by the argument in the proof of Proposition 4.6.
\end{proof}

\end{document}